\documentclass[a4p,11pt]{article}    
\usepackage{graphicx}
\usepackage[margin=2cm]{geometry}
\usepackage{amsmath,amssymb,bm}
\usepackage[T1]{fontenc}
\usepackage[utf8]{inputenc}
\usepackage{microtype}
\usepackage{newtxtext,newtxmath}
\usepackage{algorithmic,algorithm}
\usepackage{cite}
\usepackage{latexsym}
\usepackage[colorlinks=true,linkcolor=blue,citecolor=blue,urlcolor=blue]{hyperref}
\usepackage[misc,geometry]{ifsym}
\usepackage{subfigure}
\usepackage{xcolor}
\newcommand{\Eqref}[1]{Eq.~(\ref{#1})}
\newcommand{\tabref}[1]{Table~\ref{#1}} 
\newcommand{\figref}[1]{Fig.~\ref{#1}} 
\newcommand{\secref}[1]{Section~\ref{#1}}
\begin{document}

\title{Data-driven model order reduction for structures with piecewise linear nonlinearity using dynamic mode decomposition\thanks{
This manuscript is the accepted version of the article:
A. Saito and M. Tanaka,
``Data-driven model order reduction for structures with piecewise linear nonlinearity using dynamic mode decomposition,''
Nonlinear Dynamics, 111, pp. 20597--20616 (2023). 
\url{https://doi.org/10.1007/s11071-023-08958-x}
\copyright\ 2023 The Author(s). Published by Springer Nature.
}}

\author{Akira Saito\thanks{Department of Mechanical Engineering, 
Meiji University, Kawasaki, Kanagawa 214-8571, Japan. 
E-mail: asaito@meiji.ac.jp, ORCID: 0000-0001-8102-3754} 
\and 
Masato Tanaka\thanks{Toyota Research Institute of North America
    Toyota Motor North America, Inc.
    Ann Arbor, Michigan, 48105, USA, ORCID: 0000-0003-0787-0335}}

\date{}
\maketitle

\begin{abstract}
{\color{black}Piecewise-linear nonlinear systems appear in many engineering disciplines. Prediction of the dynamic behavior of such systems is of great importance from practical and theoretical viewpoint. }
In this paper, a data-driven model order reduction method for piecewise-linear systems is proposed, which is based on dynamic mode decomposition (DMD). The overview of the concept of DMD is provided, and its application to model order reduction for nonlinear systems based on Galerkin projection is explained. The proposed approach uses impulse responses of the system to obtain snapshots of the state variables. The snapshots are then used to extract the dynamic modes that are used to form the projection basis vectors. The dynamics described by the equations of motion of the original full-order system are then projected onto the subspace spanned by the basis vectors. This produces a system with much smaller number of degrees of freedom (DOFs). The proposed method is applied to two representative examples of piecewise linear systems: a cantilevered beam subjected to an elastic stop at its end, and a bonded plates assembly with partial debonding. 
The reduced order models (ROMs) of these systems are constructed by using the Galerkin projection of the equation of motion with DMD modes alone, or DMD modes with a set of classical constraint modes to be able to handle the contact nonlinearity efficiently. 
The obtained ROMs are used for the nonlinear forced response analysis of the systems under harmonic loading. It is shown that the ROMs constructed by the proposed method produce accurate forced response results. 
\end{abstract}
\section{Introduction}

{\color{black}{Piecewise-linear~(PWL) systems appear in many structural dynamics systems that are subjected to dynamic loadings. In particular, damaged structures such as with breathing cracks~\cite{SaitoEtAl2009,CasiniVestroni2011,ALShudeifatButcher2013} or delaminated composites~\cite{BurlayenkoSadowski2012} are known to be PWL systems, because they contain repetitive opening and closing of the contacting boundaries. 
Therefore, prediction of dynamic behavior of such nonlinear systems is important. 
However, the dynamics of PWL systems are known to show strong nonlinearities because of sudden changes in their stiffness, which results in rich dynamic phenomena including period doubling bifurcation and chaos~\cite{Shaw1983,Saito2018,NoguchiEtAl2022}. 
This simple yet strong nonlinearity hinders the application of standard linear analysis tools such as modal analysis. 
A general nonlinear model-order reduction (MOR) methods exist, which can reduce the size of the degrees of freedom (DOFs) of the systems~\cite{MasriEtAl1984,MasriEtAl2005}. However, its applicability to non-smooth PWL systems has never been studied.  
Therefore, there have been many attempts to develop analysis methods to predict dynamics of PWL systems in a computationally efficient manner, 
such as the ones based on nonlinear normal modes~\cite{JiangEtAl2004} and generalized bilinear amplitude approximation~\cite{TienDSouza2017,TienDSouza2019}. 
However, due to the computational cost that stems from the increasing complexity in the underlying systems with the PWL nonlinearity, application of such analysis methods is difficult. This is because of the computational cost that comes from the evaluation of the PWL nonlinearity. 
Thus, to circumvent this difficulty, the application of a \textit{data-driven} Galerkin-projection based MOR to structural dynamics problems involving PWL nonlinearity is considered in this paper. 
}}

With the increasing demand for the analysis of large-scale nonlinear systems, constructing reduced order models (ROMs) is becoming important in various engineering disciplines. That is, reduction of the number of DOFs in the system equations is of great importance to produce accurate results with reasonable amount of time, without compromising the computational accuracy. 
Therefore, many model-order reduction (MOR) methodologies have been developed. 
This paper proposes an MOR methodology based on Galerkin projection of the system equations onto a set of basis vectors obtained from {\it dynamic mode decomposition}~(DMD), which is a data-driven MOR method. 
{\color{black}{A novel approach for compensating the dynamics of modal coordinates by using classical constraint modes is also proposed to handle the contact force in the reduced order space effectively.
This combination of data-driven DMD modes and the constraint modes that are obtained based on the system model makes the proposed method a gray-box modeling approach~\cite{QuarantaEtAl2020}, instead of purely data-driven black-box modeling. 
}}

The DMD was first developed by Schmid~\cite{Schmid2010} as a method to examine coherent structures in the fluid flow field. 
DMD has then been applied to many applications involving nonlinear systems, such as acoustic mode identification in a three dimensional chamber~\cite{JourdainEtAl2013}, a swirling flow problem~\cite{BistrianEtAl2017}, combustion~\cite{RichecoeurEtAl2012}, pressure sensitive paint data~\cite{AliEtAl2016}, and many others. 
One of the greatest characteristics of the DMD over the {\it proper orthogonal decomposition~(POD)}~\cite{BerkoozEtAl1993,KerschenEtAl2005}, which is also a data-driven dimensionality reduction method, is that the DMD is able to extract not only spatial mode shapes but also their temporal information, i.e., frequency and decay-rate. If the system is linear, by using the mode shapes and the corresponding temporal information of the system for a given initial condition, DMD can be used to capture the behavior of the system by a data-driven spectral decomposition~\cite{BruntonKutz2019}, which is equation-free. 
However, if the system is nonlinear, the modes obtained by the DMD cannot be used as such an equation-free method, 
{\color {black}
{unless the DMD is applied to find a modal representation of the {\it Koopman operator} acting on carefully chosen {\it observables} of state variables such that the governing nonlinear equations becomes linear. 
Otherwise, the obtained DMD modes are linear approximations of the nonlinear manifold and the linear superposition of the DMD modes does not necessarily give accurate prediction of the original nonlinear dynamical systems. }}

Instead of equation-free approach, there have been attempts to utilize the DMD modes to construct the basis vectors for the ROMs based on {\it Galerkin projection}. 
The Galerkin projection model reduction has been successfully applied to various models such as the ones obtained by discretizing the advection diffusion equation by finite difference method~\cite{AllaKutz2017}. 
It has also been applied to a nonlinear reaction-diffusion equation~\cite{KhanNg2018}. 

In this paper, the Galerkin projection with DMD modes is applied for structural dynamics problems. There are only a few attempts to date to apply the DMD to structural dynamics problems~\cite{CunhaEtAl2022}. 
For instance, Simha et al~\cite{SimhaBiglarbegian2022} have applied the DMD to estimate natural frequencies of a linear elastic structure. 
To the knowledge of the authors, the Galerkin projection has never been applied with DMD modes for structural dynamics problems involving nonlinearity and their applicability as the basis vectors has never been discussed.

The remainder of this paper is structured as follows. In \secref{sec:method}, general projection based model order reduction procedure for equations of motion with displacement dependent nonlinearity is discussed and the theory of DMD is briefly reviewed. The DMD is then combined with the MOR procedure of the system. In \secref{section:example1}, the application of the proposed method to the forced response problem of a cantilevered beam with an elastic stop is presented. In \secref{section:example2}, the proposed method is applied to a bonded shell assembly with partial debonding, and a hybrid projection basis of DMD modes and constraint modes is examined. 
Finally in \secref{sec:conclusion}, the concluding remarks are provided. 
\section{Method}\label{sec:method}
\subsection{Model Order Reduction}
Assuming that there is an elastic body that is subject to dynamic loading and PWL nonlinearity, and that the resulting oscillation is infinitesimally small. Also assume that the governing partial differential equations are discretized by a numerical method, such as finite element~(FE) method, which results in a discrete dynamical system with piecewise-linear nonlinearity of size $m$. Denoting the nodal displacement vector as ${\bf u}(t)$, the governing equations of the system are then written as:
\begin{equation}
    {\bf M}\ddot{\bf u}(t)+{\bf C}\dot{\bf u}(t)+{\bf K}{\bf u}(t)+{\bf f}({\bf u})={\bf b}(t)\label{math:th:-2}
\end{equation}
where ${\bf M}$, ${\bf C}$, and ${\bf K}$ are mass, damping, and stiffness matrices, ${\bf f}({\bf u})$ is the nonlinear force caused by the piecewise-linear nonlinearity, ${\bf b}(t)$ is a time-dependent external forcing vector, and ${\bf M}$, ${\bf C}$, ${\bf K}\in\mathbb{R}^{m\times m}$, and ${\bf f}({\bf u})$, ${\bf b}(t)\in\mathbb{R}^m$. 
In many cases, it is desirable if the system size can be reduced by applying MOR on the system \Eqref{math:th:-2}, i.e., the Galerkin projection with a set of appropriately chosen basis vectors. 
That is, a set of vectors in $\mathbb{R}^m$, which are linearly independent, or preferably orthogonal set of vectors of dimension $m$ is selected, i.e., 
\begin{equation}
\bm{\Phi} = [\bm{\phi}_1,\dots,\bm{\phi}_p]\label{math:th:-3}
\end{equation}
where $p$ is the number of chosen basis vectors, $\bm{\phi}_i\in\mathbb{R}^m$, $i=1,\dots,p$, and in many cases, $p\ll m$. 
Then, a coordinate transformation corresponding to \Eqref{math:th:-3} is introduced, i.e., 
\begin{equation}
{\bf u}(t)\approx\sum_{i=1}^{p}\eta_i(t)\bm{\phi}_i=\bm{\Phi}\bm{\eta}(t), \label{math:th:-4}
\end{equation}
where $\eta_i(t)$ is the modal coordinate corresponding to $\bm{\phi}_i$ that acts as the variable to be solved in the reduced order space, and $\bm{\eta}(t)=[\eta_1(t),\dots,\eta_p(t)]^{\rm T}$. 
Then, the governing equations of reduced size can be obtained with respect to $\bm{\eta}(t)$, as follows: 
\begin{equation}
\tilde{\bf M}\ddot{\bm{\eta}}(t)+\tilde{\bf C}\dot{\bm{\eta}}(t)+\tilde{\bf K}\bm{\eta}(t)+\tilde{\bf f}(\bm{\eta})=\tilde{\bf b}(t), \label{math:th:-5}
\end{equation}
where $\tilde{\bf M}=\bm{\Phi}^{\rm T}{\bf M}\bm{\Phi}$, $\tilde{\bf C}=\bm{\Phi}^{\rm T}{\bf C}\bm{\Phi}$, $\tilde{\bf K}=\bm{\Phi}^{\rm T}{\bf K}\bm{\Phi}$, $\tilde{\bf b}=\bm{\Phi}^{\rm T}{\bf b}$, and $\tilde{\bf M},\tilde{\bf K},\tilde{\bf C}\in\mathbb{R}^{p\times p}$, $\tilde{\bf b}\in\mathbb{R}^p$. 
In general, the nonlinear forcing vector $\tilde{\bf f}(\bm{\eta})$ needs to be evaluated in full-order space with ${\bf u}$ especially for the PWL systems of interest, i.e., $\tilde{\bf f}(\bm{\eta})=\bm{\Phi}^{\rm T}{\bf f}(\bm{\Phi}\bm{\eta})$, because the evaluation of the piecewise linear term needs to be done in physical coordinate systems. Therefore, one needs to expand the modal coordinate vector back to full-order coordinate system at every single time step if time integration scheme is applied. This is unavoidable if one uses Galerkin projection, unless special treatment is given to the choice of the modes, as it will be shown in \secref{subsec:num2:gp}. 

It is important to choose appropriate set of $\bm{\phi}_i$'s for creating accurate ROMs especially for large-scale nonlinear problems because the nonlinear problem of interest needs to be solved in the subspace spanned by $\bm{\phi}_i$, for $i\dots p$. Therefore, the subspace should be sufficiently large such that it encompasses the dynamics of interest, yet it is small enough to conduct nonlinear calculations efficiently with reasonable accuracy. 

In this paper, the basis vectors derived from DMD are considered and the impact of the choice of $\bm{\phi}_i$'s on the resulting dynamics of Eq.~\eqref{math:th:-2} is discussed. 
Other important bases, such as linear normal modes~(LNM), which are obtained as the eigenvectors of ${\bf M}$ and ${\bf K}$ and denoted as ${\bf L}=[{\bf L}_1,\dots,{\bf L}_p]$, ${\bf L}_i\in\mathbb{R}^m$, and eigenvectors obtained from the POD, which will be discussed later, are also taken into account for comparison. 
\subsection{Dynamic Mode Decomposition}
In this section, mathematical background of the DMD is briefly reviewed. 
First, let us assume that ${\bf x}(t)\in\mathbb{R}^m$ denotes a time-dependent vector of generalized coordinates, or a state vector of an $m$-dimensional dynamical system. Now assume that a snapshots of the generalized coordinates are taken at discrete time instants $t_j$ for $j=1,\dots,n$, and are stored in a single matrix so called {\it snapshot matrix}, ${\bf X}=[{\bf x}_1,\dots,{\bf x}_n]$ for ${\bf x}_j\triangleq{\bf x}(t_j)$ and $t_{j+1}=t_{j}+\Delta t$ where the size of $\Delta t$ is fixed. 
Next, in DMD formulation, another snapshot matrix ${\bf Y}$ is considered, where the columns of ${\bf Y}$ are the time-shifted versions of ${\bf X}$, i.e., ${\bf Y}\triangleq[{\bf x}_2,\dots,{\bf x}_{n+1}]$. 
The key assumption we make here is that there exists a linear map represented by a matrix ${\bf A}$ between ${\bf x}_j$ and ${\bf x}_{j+1}$, i.e., 
\begin{equation}
	{\bf x}_{j+1}={\bf A}{\bf x}_{j}\label{eq:th:-1}. 
\end{equation}
This yields, 
\begin{equation}
	{\bf Y} = {\bf A}{\bf X}. 
\end{equation}
It means that the following relationship holds: 
\begin{equation}
{\bf A}={\bf Y}{\bf X}^{\dagger}\label{eq:th:-1+1}, 
\end{equation}
and ${\bf X}^{\dagger}$ denotes the Moore-Penrose pseudo inverse matrix of ${\bf X}$. 
The {DMD modes and eigenvalues} are defined as the eigenvectors and eigenvalues of ${\bf A}$\cite{Schmid2010}. 
By applying singular value decomposition~(SVD) to {\bf X}, ${\bf X}={\bf U}\bm{\Sigma}{\bf V}^*$ where ${\bf U}\in\mathbb{C}^{m\times r}$,  $\bm{\Sigma}\in\mathbb{R}^{r\times r}$, and ${\bf V}\in\mathbb{C}^{n\times r}$, $r$ is the rank of ${\bf X}$, and $^*$ denotes the Hermitian transpose. 
Note that ${\bf U}$ contains the eigenvectors of the POD~\cite{KerschenEtAl2005}, or the {\it POD modes}, and they are designated as $[{\bf U}_1,\dots,{\bf U}_{\color{black}{r}}]$. 
The diagonal terms in $\bm{\Sigma}$ are denoted as $\sigma_i$, $i=1,\dots,r$.
The pseudo inverse ${\bf X}^\dagger$ can then be computed by using ${\bf V}$, ${\bf U}$, and $\bm{\Sigma}$, as follows: 
\begin{equation}
{\bf X}^{\dagger}={\bf V}\bm{\Sigma}^{-1}{\bf U}^*. 
\end{equation}
The DMD eigenvalue is now obtained by computing the eigenvalues of the projected matrix ${\bf A}$ onto the POD modes~\cite{TuEtAl2014}. Namely, defining the projected matrix as 
\begin{align}
\tilde{\bf A}&\triangleq{\bf U}^*{\bf A}{\bf U},\nonumber\\
&={\bf U}^*{\bf Y}{\bf X}^{\dagger}{\bf U},\nonumber \\
&={\bf U}^*{\bf Y}{\bf V}\bm{\Sigma}^{-1}{\bf U}^*{\bf U},\nonumber \\
&={\bf U}^*{\bf Y}{\bf V}\bm{\Sigma}^{-1}. 
\end{align}
Now the DMD eigenvalues of $\tilde{\bf A}$ can be obtained by solving the following eigenvalue problem, i.e.,
\begin{equation}
\tilde{\bf A}{\bf w}_i=\mu_i{\bf w}_i, 
\end{equation}
where $\mu_i$  is the $i$th DMD eigenvalue. 
The DMD mode $\bm{\varphi}_i$ corresponding to $\mu_i$ of dimension $m$ is then defined by, as in the definition of the {\it exact DMD} in Ref.~\cite{TuEtAl2014}:, 
\begin{equation}
\bm{\varphi}_i = ({\bf Y}{\bf V}\bm{\Sigma}^{-1}){\bf w}_i.
\label{eq:th:32}
\end{equation}
This holds because $\bm{\varphi}_i$ is the eigenvector of ${\bf A}$ corresponding to $\mu_i$, i.e., 
\begin{align}
{\bf A}\bm{\varphi}_i&={\bf Y}{\bf V}\bm{\Sigma}^{-1}{\bf U}^*{\bf Y}{\bf V}\bm{\Sigma}^{-1}{\bf w}_i,\nonumber\\
&={\bf Y}{\bf V}\bm{\Sigma}^{-1}\tilde{\bf A}{\bf w}_i,\nonumber\\
&={\bf Y}{\bf V}\bm{\Sigma}^{-1}\mu_i{\bf w}_i, \nonumber\\
&=\mu_i\bm{\varphi}_i. 
\end{align}
It is noted that most DMD eigenvalues appear with their complex-conjugate pairs, so do the corresponding DMD modes. Therefore, both real and imaginary parts are used for constructing the projection basis. 
%
It is also known that the DMD yields an approximate eigen-decomposition of the best-fit linear operator relating two data matrices ${\bf X}$ and ${\bf Y}$~\cite{Tu2013}, which in this case is the matrix ${\bf A}$. 
One of the advantages of the DMD over POD is that the DMD eigenvalues can be associated with their frequency and decay rate, or damping ratios. Namely, defining $s_i$, $i=1,\dots r$, as 
\begin{equation}
    s_i =\log(\mu_i)/\Delta t,
\end{equation}
which is also written as $s_i=-\zeta_i\omega_i{\color{black}{\pm}}{\rm j}\omega_i\sqrt{1-\zeta_i^2}$ where ${\rm j}=\sqrt{-1}$ {\color{black}{with the assumption that all DMD modes decay and $0\leqslant\zeta_i\leqslant 1$}}. Then, the following relationships hold.  
\begin{align}
    f_i&=\omega_i/2\pi=|s_i|/2\pi,\\
    \zeta_i&
    =-{\rm Re}({s_i})/|s_i|,
\end{align}
where {\color{black}{$\omega_i$}}, $f_i$ {\color{black}{and $\zeta_i$ are}} the {\color{black}{angular frequency}}, frequency {\color{black}{and damping ratio }}corresponding to the {\color{black}{$i$-th}} DMD mode{\color{black}{, respectively.}}
Although they are error prone, it was shown that $f_i$ and $\zeta_i$ coincide with the undamped natural frequency and modal damping ratio if the matrices ${\bf X}$ and ${\bf Y}$ come from the response of a linear system~\cite{SaitoKuno2020}. If the system is nonlinear, they are not necessarily the natural frequencies and modal damping ratios. 
Furthermore, the magnitude of $\mu_i$ can be used to rank the DMD modes based on its dominance in the snapshot matrices. Namely, 
\begin{align}
\mu_i
&={\rm exp}{\left\{{\rm Re}(s_i)\Delta t{\color{black}{\pm}}{\rm jIm}(s_i)\Delta t\right\}},\nonumber \\
&={\rm exp}\left\{-\zeta_i\omega_i\Delta t\right\}
  {\rm exp}\left\{{\color{black}{\pm}}{\rm j}\omega_i\sqrt{1-\zeta_i^2}\Delta t\right\}.
\end{align}
Therefore, 
\begin{equation}
|\mu_i|=
|{\rm e}^{-\zeta_i\omega_i\Delta t}|
=1/|{\rm e}^{\zeta_i\omega_i\Delta t}|.
\end{equation}
This means that if $|\mu_i|$ is large, then $\zeta_i\omega_i$ is small for a fixed $\Delta t$, i.e., the decay rate of the oscillation of the corresponding DMD mode component is small. In other words, it decays slowly. If $|\mu_i|$ is small, then $\zeta_i\omega_i$ is large, i.e., the decay rate of the oscillation of the corresponding component is large. It means that it decays fast. Therefore, when constructing a ROM using $\bm{\varphi}_i$ by the Galerkin projection, it is desirable to include modes with low decay rate, or, with large $|\mu_i|$ as many as necessary, because they are expected to dominate the dynamics. 

%
\subsection{Analysis Procedure}
In this paper, forced response of structures with PWL nonlinearity is considered. 
The process of the analysis based on DMD-based MOR is threefold: (1) computation of the snapshots, (2) ROM construction with the DMD modes, and (3) the nonlinear forced response calculation with the ROM. 

First, in order to obtain the DMD-based ROMs, snapshots $[{\bf x}(t_1),{\bf x}(t_2),\dots,{\bf x}(t_{n+1})]$ need to be obtained numerically or experimentally. Suppose that it is obtained numerically, the difficulty of the procedure of the computation of the snapshots is problem dependent. Also, if the snapshot is obtained by solving initial value problems, it is important to set initial conditions or initial loading such that they excite sufficiently many DMD modes that are possibly involved in the dynamics to be predicted with the ROM. 

Second, with the snapshot matrices, DMD needs to be conducted. The challenge here is that the DMD spectrum may contain DMD eigenvalues that correspond to fast-decaying modes or non-existing modes that do not represent any physical dynamic phenomena. Therefore, such DMD modes need to be removed. In this paper, DMD mode selection is conducted based on singular value rejection~\cite{ThiteThompson2003}, and successive application of DMD with increasing sampling frequency~\cite{SaitoKuno2020}. 
Note that the sampling frequency of the original time history cannot be changed once it is computed. Therefore, it is re-sampled with larger sampling periods, and the DMD is conducted to see whether the DMD eigenvalues consistently appear regardless of the sampling frequency. This is inspired by the concept of the stability diagram~\cite{PeetersEtAl2004} that is used in the frequency-domain pole estimator used in the area of experimental modal analysis. 
Moreover, the computed eigenvalues are ranked based on their magnitudes to further promote the sparcity of the DMD eigenvalues, which results in the reduction of the system size in the ROM. 

Third, the obtained DMD mode shapes $\bm{\varphi}_i$'s are used to form a transformation matrix $\bm{\Phi}=[\bm{\varphi}_1,\dots,\bm{\varphi}_p]$. The reduced system matrices and reduced forcing vectors are then obtained. 
There are possibly two methods to handle nonlinear force that comes from the PWL nonlinearity. The first approach is to compute the nonlinear force in full-order space and project it back to reduced space. Namely, it is schematically written as follows. 
\begin{equation}
{\bf f}={\bf f}(\bm{\Phi}\bm{\eta}(t))\rightarrow\tilde{\bf f}=\bm{\Phi}^{\rm T}{\bf f}.\label{math:eq21}
\end{equation}
This process requires the expansion of the modal coordinates $\bm{\eta}(t)$ back to larger set of physical DOFs ${\bf u}(t)$ at every time step, because all modes in $\bm{\Phi}$ contribute to the nonlinear force. The second approach is to compute the nonlinear force in reduced-order space directly, without expanding the modal coordinates back to physical DOFs. In this paper, it is shown that this can be achieved by adding {\it constraint modes}\cite{Craig1968} in $\bm{\Phi}$ so that the DOFs subjected to the nonlinearity are accessible in the reduced-order subspace. The process is schematically written as follows. 
\begin{equation}
\tilde{\bf f}=\tilde{\bf f}(\bm{\eta}(t)). \label{math:eq22}
\end{equation}
This process is computationally more efficient than the first approach represented as Eq.~\eqref{math:eq21}, because the multiplications of $\bm{\Phi}^{\rm T}\bm{\eta}$ and $\bm{\Phi}^{\rm T}{\bf f}$ are not necessary. 
The overall MOR algorithm explained above is shown in Algorithm~\ref{alg1}. 

\begin{figure*}[tb]
\centering
\includegraphics[scale=1]{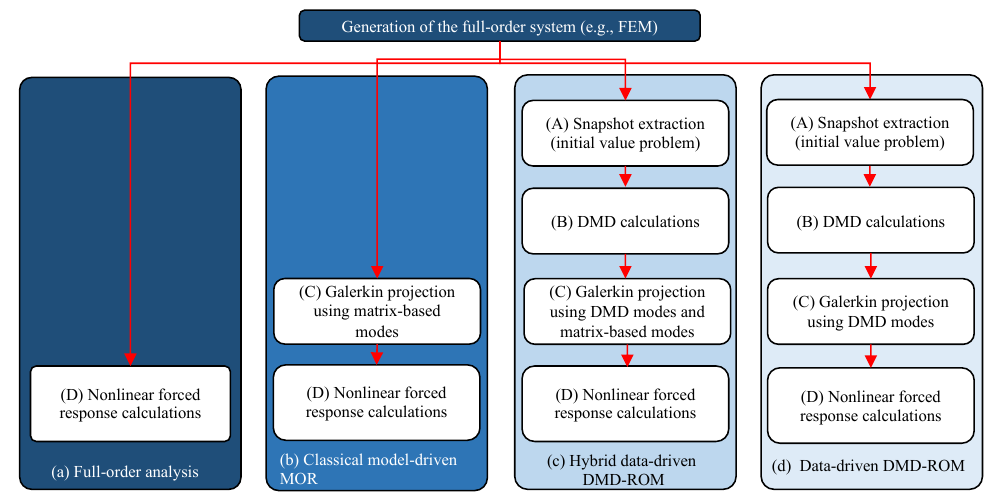}
\caption{Overview of the forced response analysis}
\label{fig:overview}
\end{figure*}
{\color{black}To better visualize the structure of the proposed methodology, analysis flowcharts of typical forced response calculations are shown in \figref{fig:overview}. 
First, one needs to generate the full-order system by discretizing the spatial variables of the governing equations by numerical methods, such as FEM. A naive approach for obtaining the forced response solution is to run the forced response with full-order system, as shown in \figref{fig:overview}(a). One may want to apply classical model-based MOR to the system that are based on superposition of modes that can be obtained by modal and static analyses of the system matrices with special boundary conditions, such as Craig-Bampton matrix condensation~\cite{Craig1968}, as shown in \figref{fig:overview}(b). This process consists of the Galerkin projection of the equation onto the set of modes, and the nonlinear forced response calculations of the reduced system. This has been a good approach for gaining computational efficiency without degrading the accuracy. However, the modes obtained purely from static or modal analyses are from linear analyses. Therefore, they may not necessarily span the subspace where the nonlinear response lies. 
Data-driven MOR, which was proposed in Ref.~\cite{AllaKutz2017} and is extended to PWL system in this paper, on the other hand, is expected to be able to use a good basis vectors for nonlinear forced response calculations, because they come from DMD of nonlinear response of an initial value problem of the system to be studied in the forced response calculations. This process is shown in \figref{fig:overview}, which consists of the extraction of snapshots and DMD calculations from time series data. The Galerkin projection is then applied to the sytem equations to obtain the ROM. Nonlinear forced response is then computed. This is considered in \secref{section:example1}. 
Furthermore, when one forms ROMs for PWL systems, handling the contact force at the contacting boundaries can be expensive even with ROMs either by classical model-based MOR or data-driven MOR, because the evaluation of the contact force needs to be done in physical domain, which results in repetitive computation of modal coordinates back to physical coordinates at every time instant. Given the number of DOFs that are subject to PWL nonlinearity is not prohibitively large, one is able to efficiently compute the nonlinear force in the reduced subspace by keeping the physical DOFs in the ROM, by using static modes called constraint modes. This process is shown in \figref{fig:overview}(c), which is denoted here as Hybrid data-driven MOR, and this approach is considered in \secref{section:example2}. 
}

{\color {black}{
There have been a couple of methods to apply projection-based methods to nonlinear systems. 
For instance, Discrete Empirical Interpolation Method~(DEIM), which is based on POD, has been proposed to reduce the number of DOFs where nonlinear terms are evaluated~\cite{ChaturantabutEtAl2010}. 
The method has been extended to be used with DMD as well~\cite{AllaKutz2017,HesthavenEtAl2022}. 
Hyper-reduction has also been proposed for enhancing the numerical efficiency of the ROM by using a set of sampling grid points in the computational domain to specify reduced integration domains where the problem is solved~\cite{FritzenEtAl2018}.  
The problems studied in this paper are assumed to be fundamentally {\it localized}, i.e., the size of the region where the PWL nonlinearity appears is relatively small in comparison with the entire computational domain, yet it is subjected to the strong PWL nonlinearity. Hence, these methods that are specifically tailored for reducing the number of DOFs that are subject to nonlinearities do not have to be applied to the systems of interest in this study. 
}}

The following sections provide numerical examples of the application of the proposed approach to the forced response problems of mechanical systems that are subject to piecewise-linear nonlinearities. 

\begin{algorithm}[bt]
\begin{algorithmic}[1]
\caption{Model order reduction based on Dynamic Mode Decomposition}\label{alg1}
\STATE{Obtain $[{\bf x}_1,{\bf x}_2,\dots,{\bf x}_{n+1}$] for $t_{j+1}=t_{j}+\Delta t$} by solving initial value problem for a given initial condition
\FOR{$k=1,\dots k_{max}$}
\STATE{$\Delta t_k\leftarrow (k_{max}/k)\Delta t$}
\STATE{Re-sample the snapshots with $\Delta t_k$ and obtain $[{\bf x}_1,{\bf x}_2,\dots,{\bf x}_{n+1}$] for $t_{j+1}=t_{j}+\Delta t_k$}
\STATE{${\bf X}\leftarrow[{\bf x}_1,{\bf x}_2,\dots,{\bf x}_n]$ and ${\bf Y}\leftarrow[{\bf x}_2,{\bf x}_3,\dots,{\bf x}_{n+1}]$}
\STATE {\color{black}{Compute SVD: ${\bf U}\bm{\Sigma}{\bf V}^*\leftarrow{\bf X}$}}
\STATE{$\tilde{\bf A}\leftarrow{\bf U}^{*}{\bf Y}{\bf V}\bm{\Sigma}^{-1}$}
\STATE{Solve $\tilde{\bf A}{\bf w}_i=\mu_i{\bf w}_i$}
\STATE{$s_i\leftarrow\mathrm{log}\left(\mu_i\right)/\Delta t_k$, $f_i\leftarrow|s_i|/2\pi$}
\STATE{$\bm{\varphi}_i=\left({\bf Y}{\bf V}\bm{\Sigma}^{-1}\right){\bf w}_i$}
\ENDFOR{}
\STATE{Form $\bm{\Phi}=[\bm{\varphi}_1,\dots,\bm{\varphi}_p]$}
\STATE{Compute 
$\tilde{\bf M}=\bm{\Phi}^{\rm T}{\bf M}\bm{\Phi}$, 
$\tilde{\bf K}=\bm{\Phi}^{\rm T}{\bf K}\bm{\Phi}$, 
$\tilde{\bf C}=\bm{\Phi}^{\rm T}{\bf C}\bm{\Phi}$, 
$\tilde{\bf b}=\bm{\Phi}^{\rm T}{\bf b}(t)$}
\end{algorithmic}
\end{algorithm}
\section{Forced response of a beam with an elastic stop
{\color{black}{with data-driven model-order reduction}}
}\label{section:example1}
\begin{figure}[tb]
    \centering
    \includegraphics[scale=1]{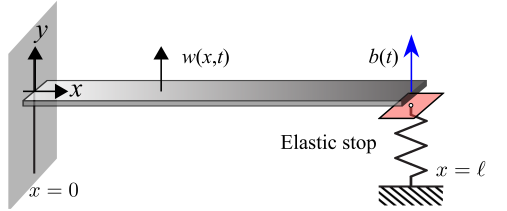}
    \caption{Cantilevered beam subject to harmonic forcing and elastic stop}\label{numerical_examples:fig1}
\end{figure}
The approach is applied to a cantilevered beam with an elastic stop at its end, and is subject to an external forcing. The schematics of the problem is shown in \figref{numerical_examples:fig1}. The governing equation of motion of the beam is obtained with the Euler-Bernoulli beam theory, in strong form as follows, 
\begin{equation}
    \rho A\frac{\partial^2 w}{\partial t^2}+
    EI\frac{\partial^4 w}{\partial x^4}=0,
\end{equation}
where $E$ is the Young's modulus, $I$ is the second moment of area, $\rho$ is the density, $A$ is the cross-sectional area of the beam. 
The boundary conditions at $x=0$ and $x=\ell$ are, 
\begin{align}    
    w(0,t)&=0,\label{math:eq:bc1}\\
    EI\frac{\partial^3 w}{\partial x^3}(\ell,t)&=
    b(t)+f_{nl}(w(\ell,t)),
\end{align}
where $b(t)$ is the external harmonic forcing, and $f_{nl}(w(\ell,t))$ is the PWL nonlinear restoring force acted upon by the elastic stop. 
Namely, the spring generates the restoring force only when $w(\ell,t)\leqslant0$, i.e., 
\begin{equation}
f_{nl}(w(\ell,t)) = -k_c{\rm max}\left(
-w(\ell,t),0
\right), \label{eq:beam_nl_force}
\end{equation}
where ${\rm max}(x,y)$ returns $x$ if $x> y$, $y$ if $x< y$, or $x(=y)$ if $x=y$, and $k_c$ is the spring constant of the elastic stop. 
With these governing equations, the governing equation was discretized by Euler-Bernoulli beam elements, which resulted in 32 elements of 64~DOFs where each node has a transverse and a rotational DOFs. 
The fixed boundary condition of \Eqref{math:eq:bc1} has been achieved by removing the DOFs at the root. 
The resulting equations of motion have the form: 
\begin{equation}
{\bf M}\ddot{\bf w}(t)+{\bf C}\dot{\bf w}(t)+{\bf K}{\bf w}(t)={\bf b}(t)+{\bf f}({\bf w})\label{num_gov_eq1}
\end{equation}
where ${\bf w}$ is the displacement vector containing all transverse and rotational DOFs of all nodes. 
\subsection{Snapshot extraction and dynamic mode decomposition}
\begin{figure}[tb]
\centering
\subfigure[Impulsive force applied at the tip]{\includegraphics[width=8.5cm]{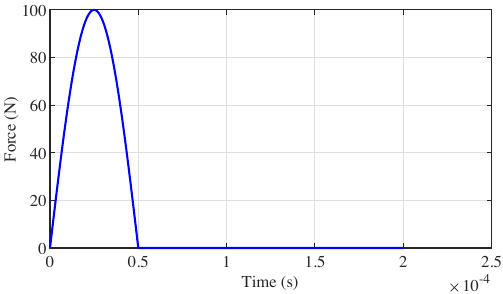}}
\subfigure[Response of the beam at the tip]{\includegraphics[width=8.5cm]{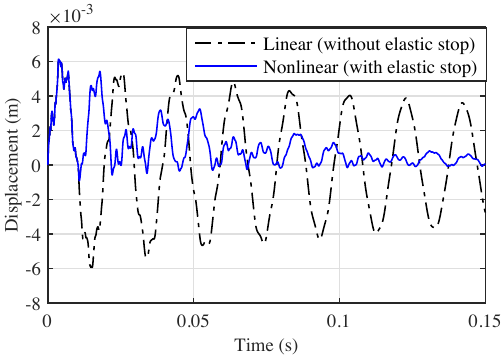}}
\caption{Applied force and the response of the beam for both linear and nonlinear cases}
\label{numerical_example:fig1}
\end{figure}

As discussed in \secref{sec:method}, snapshots should be obtained to compute DMD or even POD, which is the characteristics of these data-driven approach. The choice of snapshots is arbitrary, but they need to contain rich information about the dynamics to be captured by the ROM using DMD modes. 
Therefore, for this numerical example, snapshots were obtained by numerically solving \Eqref{num_gov_eq1} for an impulsive force, i.e., {\it impulse response}, which is expected to excite the motion of the beam for a wide frequency range. 
The impulse response of the beam can be obtained by applying impulsive force at the tip of the beam and the time response of the system was computed by a time integration method. 
The applied impulsive force was defined as a sinusoidal wave of half-period define as, 
\begin{equation}
b(t)=
\left\{
\begin{array}{ll}
b_0\sin\left(
2\pi t/T\right),& 0\leqslant t \leqslant T/2\\
0,& t> T/2
\end{array}
\right.
\end{equation}
where $T$ is the period of the sinusoid with $T=1.0\times10^{-4}$s and $b_0=100$N. The value of $T$ was chosen such that it is smaller than the period of the fastest dynamics to be considered in this problem, which in this case was judiciously chosen to be 4000~Hz in frequency or $2.5\times10^{-4}$s in time duration. 
The time history of the applied force is shown in \figref{numerical_example:fig1}(a). 
The equation of motion \Eqref{num_gov_eq1} with ${\bf w}(0)={\bf 0}$ and $\dot{\bf w}(0)={\bf 0}$ was solved by \verb?ode45? solver implemented in Matlab\textregistered. The value of the spring constant for the elastic stop was chosen to be $k_c=1000$N/m. 
The obtained impulse response of the beam subject to the elastic stop measured at the tip is shown in \figref{numerical_example:fig1}(b). 
To see the effects of the elastic stop at the tip, the response without the elastic stop, or the case where $k_c=0$ is labeled as Linear in \figref{numerical_example:fig1}(b). 
As can be seen, both responses are identical until the tip strikes the elastic stop where $w=0$. After the first strike, the tip displacement changes its direction due to the elastic stop, which cannot be seen in the linear response. This process is repeated until this bouncing behavior stops when the kinetic energy disappears. 
It is noted that higher frequency components are excited in the nonlinear response than in the linear response. 
\begin{figure}[tb]
\subfigure[Linear (no elastic stop)]{\includegraphics{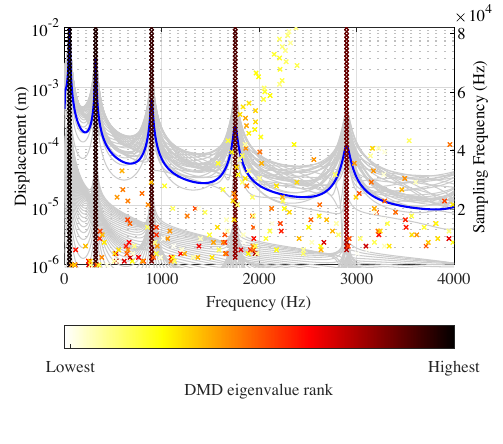}}
\subfigure[Nonlinear (with elastic stop)]{\includegraphics{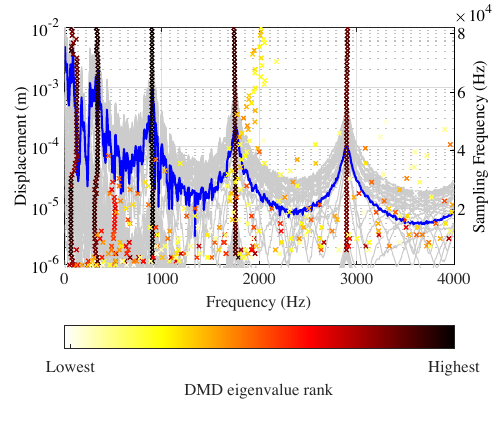}}
\caption{Pseudo-stability diagram of the DMD eigenvalues, along with the FFT results of the displacments. {\color{lightgray}{\bf ---}}: FFT spectra of all DOFs, {\color{blue}{\bf ---}}: average spectrum, ${\bf \times}$: DMD eigenvalues.}\label{numerical_example:fig2}
\end{figure}

Next, using the obtained snapshots of the impulse response, DMD has been applied. Based on the algorithm shown in Algorithm.~\ref{alg1}, DMD modes were extracted for both linear and nonlinear cases. The results are shown as the pseudo-stability diagram of the DMD~\cite{SaitoKuno2020}. 
In the pseudo-stability diagram, DMD eigenvalues with increasing sampling frequencies and the spectra of the displacement computed by the Fast Fourier Transform~(FFT). The color of the DMD eigenvalue represents the ranking of the eigenvalue among the ones computed for a specific sampling frequency. The larger the magnitude of the DMD eigenvalue is, the higher its ranking is. 
From this diagram, we can see if a DMD eigenvalue and the corresponding DMD mode shape should be taken into account in the ROM, i.e., if the DMD eigenvalue appears with high ranking regardless of the sampling frequency, then the eigenvalue is stable and hence is important. It means that it is likely that it needs to be taken into account in the ROM.  
\begin{table}[tb]
\centering
\caption{DMD eigenvalues in Hz for linear and nonlinear cases}\label{numerical_example:tab1}
\begin{tabular}{cccc}\hline
Mode number & Linear & Nonlinear & Deviation (\%)\\ \hline
1 & 50.9 &65.3 &   28.29\\
2 & 320 &316 &   -1.250\\
3 & 893&886&    -0.7839\\
4 & 1745&1735 &   -0.5731\\
5 & 2886&2874 &   -0.4158\\ \hline
\end{tabular}
\end{table}
\begin{figure}[tb]
\centering
\subfigure[LNM]{\includegraphics[]{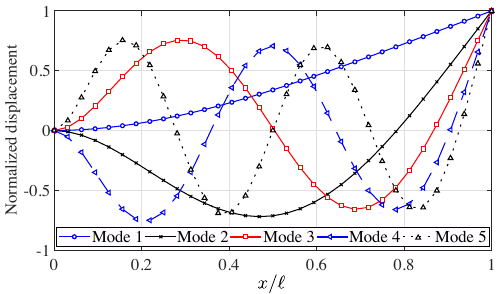}}
\subfigure[POD modes]{\includegraphics[]{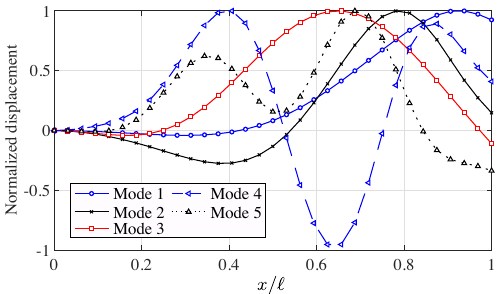}}
\subfigure[DMD modes]{\includegraphics[]{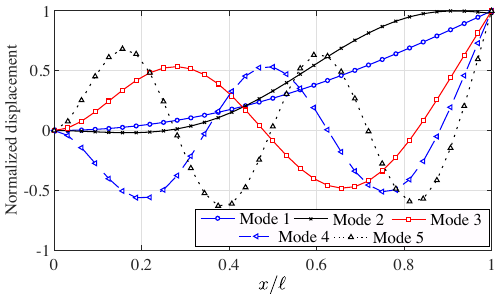}}
\caption{Comparison between LNM, POD modes and DMD modes}
\label{fig:lnm_vs_dmd}
\end{figure}

Based on the diagrams shown in \figref{numerical_example:fig2}, we can see that there are five dominant DMD eigenvalues for both linear and nonlinear cases. The DMD eigenvalues correspond to the peak frequencies of the FFT spectra, especially for the linear case. 
In fact, it is known that DMD is equivalent to the Ibrahim's time domain modal parameter extraction method~\cite{Ewins2009,SaitoKuno2020}. Therefore, if the DMD is applied to the linear response, 
{\color{black}{the obtained DMD modes match exactly the linear normal modes obtained from the mass and stiffness matrices, and the corresponding DMD eigenvalues contain the undamped natural frequencies. Furthermore, the modal damping ratio obtained from the DMD eigenvalues exactly match the modal damping ratios obtained from the mass, stiffness, and damping matrices of the system if the damping is modeled as proportional damping where the damping term can also be diagonalized by the linear normal modes.}}
On the other hand, as can be seen in \figref{numerical_example:fig2}(b) for the nonlinear case, the FFT spectra contain many frequency components and clear dominant peaks cannot be determined from the FFT alone, especially for the low frequency range. However, the DMD eigenvalues coincide with the peak frequencies of the FFT spectra. 
The obtained frequencies corresponding to the DMD eigenvalues are shown in \tabref{numerical_example:tab1}. Note that the frequencies extracted from the linear response are identical to the natural frequencies that can be obtained from the mass and stiffness matrices of the system. On the contrary, the DMD eigenvalues extracted from the nonlinear response do not match the natural frequencies as expected. The most significant shift in the eigenvalue is 28.29\% for the first mode in comparison with that of the linear case. This appears to be the stiffening effect due to the elastic stop at the tip. The deviations for all the other modes are much smaller than that for the first mode. 

\begin{figure}[tb]
\centering
\includegraphics{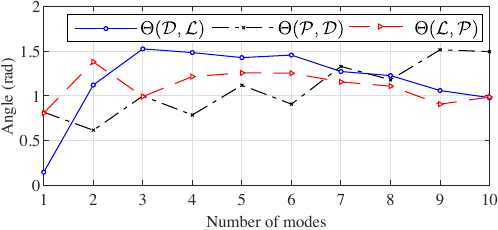}
\caption{Principal angles between the subspaces spanned by LNMs, POD modes, and DMD modes for different number of modes}\label{numerical_example:fig31}
\end{figure}
Next, the obtained mode shapes are discussed. First five modes are shown in \figref{fig:lnm_vs_dmd} for LNMs, POD modes, and the DMD modes. 
For the POD and DMD modes, snapshots obtained from the nonlinear case were used where the elastic stop at the tip exists. 
Each mode shape is normalized with respect to its infinity norm. 
As can be seen from \figref{fig:lnm_vs_dmd}(a) and (c), despite the slight differences in their amplitude values, the DMD modes resemble the LNMs for the first five modes, except the second mode. 
On the other hand, the POD modes do not resemble neither LNMs nor DMD modes for the first five modes, as seen in \figref{fig:lnm_vs_dmd}(b). 

Next, to see the nature of the subspaces spanned by the modes, the principal angles ~\cite{GolubLoan1996} between the subspaces are examined. Let $\mathcal{L}$, $\mathcal{P}$, and $\mathcal{D}$ denote the subspaces spanned by the set of LNMs, POD modes, and DMD modes, respectively, or 
$\mathcal{L}={\rm Span}({\bf L}_1,\dots,{\bf L}_p)$, 
$\mathcal{P}={\rm Span}({\bf U}_1,\dots,{\bf U}_p)$, 
$\mathcal{D}={\rm Span}(\bm{\varphi}_1,\dots,\bm{\varphi}_p)$. 
Then, the principal angles between these subspaces were computed and shown in \figref{numerical_example:fig31}, where $\Theta(\mathcal{X},\mathcal{Y})$ designates the principal angle between two subspaces $\mathcal{X}$ and $\mathcal{Y}$. 
From \figref{numerical_example:fig31}, we can see that the principal angle between $\mathcal{D}$ and $\mathcal{L}$ is small when the number of modes is one, which means they are almost linearly dependent. This agrees with the observation that the mode shapes of the first mode for both LNM and DMD modes were almost identical to each other. 
However, as the number of modes increases, the principal angle reaches almost $\pi/2$ when the number of modes is three, which means that they become almost orthogonal to each other, then slightly decreases down to 1.0. 
The pricipal angle between $\mathcal{P}$ and $\mathcal{D}$ gradually increases and become close to $\pi/2$ when the number of modes is nine. It means that they are almost orthogonal with each other. 
On the other hand, the principal angle between $\mathcal{P}$ and $\mathcal{L}$ is between approximately 0.5 and 1.5, which means their relationship lies between linearly dependent and orthogonal. 

From these results, we can see that the subspaces spanned by these sets of vectors are not parallel to each other, i.e., LNMs, PODs, and DMDs do not form the same subspace. Therefore, their capability to express the nonlinear dynamics should also be different from each other, and needs to be further examined with forced response calculations, as follows. 

\subsection{Forced response calculation using reduced order models}
\begin{figure*}[tb]
    \centering
\includegraphics[scale=1]{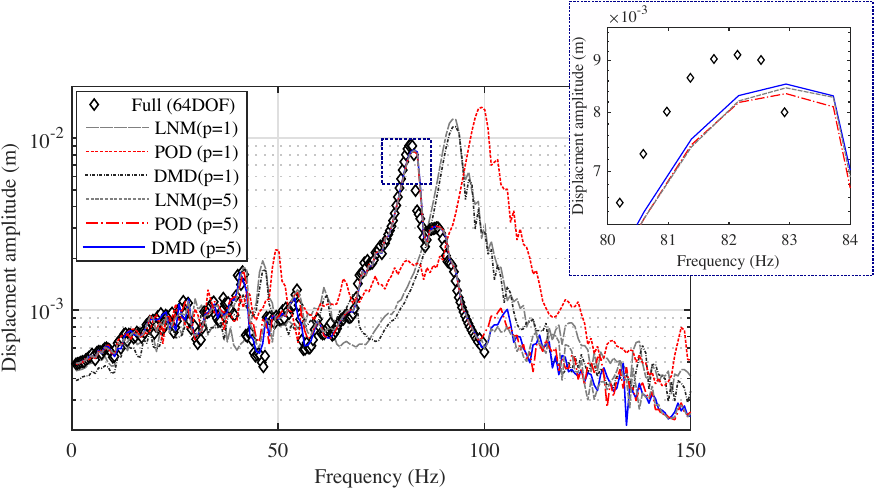}
    \caption{Forced response results of full-order FEM, POD, and DMD}\label{numerical_example:fig4}
\end{figure*}
Forced response calculations were conducted using full FE models, reduced order models with LNMs, POD modes, and DMD modes, and the results are examined. 
For the reduced order models with DMD modes, the MOR procedure follows the Algorithm~\ref{alg1}. For the other modes, the process of the Galerkin projection is rather straightforward, i.e., with $\bm{\Phi}={\bf U}=[{\bf U}_1,\dots,{\bf U}_p]$ for POD modes, or $\bm{\Phi}={\bf L}=[{\bf L}_1,\dots,{\bf L}_p]$ for LNMs and $\tilde{\bf M}=\bm{\Phi}^{\rm T}{\bf M}\bm{\Phi}$, 
$\tilde{\bf K}=\bm{\Phi}^{\rm T}{\bf K}\bm{\Phi}$, 
$\tilde{\bf C}=\bm{\Phi}^{\rm T}{\bf C}\bm{\Phi}$. 
A harmonic forcing was applied at the beam tip, i.e.,
\begin{equation}
    b(t)=A\sin(\omega t)
\end{equation}
where $A=0.01$N, and the steady-state response of the system was sought for the frequency range of $\omega/(2\pi)=0$ to 120~Hz to excite a resonance that was expected to occur near the first linear natural frequency of 50.9~Hz. 
The corresponding external forcing vector is evaluated as $\tilde{\bf b}=\bm{\Phi}^{\rm T}{\bf b}(t)$. The nonlinear force corresponding to Eq.~\eqref{eq:beam_nl_force} is evaluated with physical DOFs that are transformed back from the modal coordinates. 
\begin{figure}[tb]
\centering
\includegraphics[scale=1]{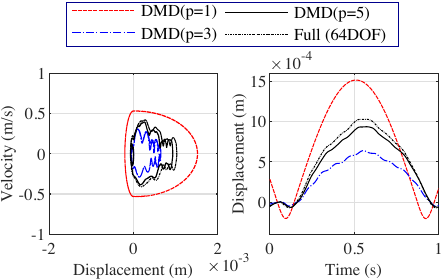}
\caption{Phase portrait of the tip movement at the resonance}
\label{numerical_example:fig6}
\end{figure}
The results of the forced response calculations with full-order FEM model, the ROMs constructed by LNMs, POD modes, and DMD modes for $p=1$ and $5$ are shown in \figref{numerical_example:fig4}. 
The displacement was evaluated at the tip of the beam. 
As can be seen, there is a resonant peak at 82~Hz for the FEM model. All ROMs with $p=1$ do not capture neither the maximum amplitude value nor the frequency of the resonant frequency, as shown in \figref{numerical_example:fig4}. The ROMs with the LNM and the DMD mode better capture the resonance than the ROM with the POD modes, but they are still lower than the one predicted by the full model by 10\%. 
On the other hand, the ROMs with $p=5$ capture this resonance behavior quite well for all ROM types. This means that all three projection basis vectors are equivalently good to capture the resonant behavior if sufficiently large number of modes are kept in the ROM. In particular, the resonance predicted by the DMD modes is the best among the ones considered, as seen in the enlarged plot shown in \figref{numerical_example:fig4}. 
{\color{black}{The results of full order and the ROM are considered to match quantitatively in this case, considering that the errors between the nonlinear resonant frequencies obtained by the full order system and the ROM is approximately 1\%, and the errors in the peak amplitude are approximately 6\%. With these, we belive that the accuracy of the ROM is satisfactory.}}
%

%
Moreover, to see the behaviors of the systems at the resonance, the time histories of the tip displacement and the corresponding phase portraits are examined for the ROMs with the DMD alone. 
Figure \ref{numerical_example:fig6} shows the time history of the tip displacement at the resonance frequency. 
As expected, ROMs with $p=1$ and $p=3$ do not capture the trajectory of the full-order system well. 
On the other hand, the ROM with $p=5$ captures the temporal behavior of the tip displacement computed by the full-order FEM model. 

From these results shown above, we can see that DMD is able to capture the dynamics of a simple beam structure that is subjected to an elastic stop, which causes the system to have piecewise-linearity. The accuracy of the ROM is slightly better than those of LNM and POD modes, which have already been known to be good basis vectors for MOR of linear and nonlinear structural dynamics problems. 
\section{Forced response of a bonded shell assembly with partial debonding {\color{black}{with hybrid data-driven model-order reduction}}}\label{section:example2}
\begin{figure}[tb]
\centering
\includegraphics[width=8.5cm]{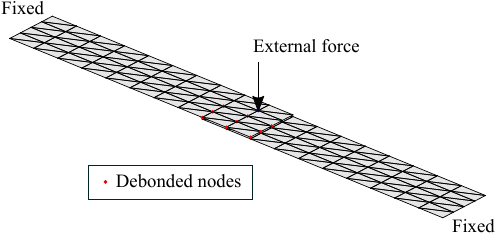}
\caption{Bonded shell assembly with partial debonding}
\label{numerical_example2:fig1}
\end{figure}
This section provides a numerical example of the forced response problem of a bonded shell assembly with partial debonding. 
The FE model of the assembly is shown in \figref{numerical_example2:fig1}. 
Namely, two plates are bonded together at one of their ends, but 75\% of the bonded area is assumed to be debonded. This causes the plates to contact with each other repeatedly if they are subject to periodic loading, such as harmonic forcing, which results in nonlinear dynamic behavior. 
The shell assembly has been discretized by MITC3 isotropic triangular shell elements~\cite{Lee2004}, where each node of the element has three translational and two rotational DOFs. The total number of DOFs of the system is 500. 
Both ends are assumed to be fixed, and an external force is assumed to be applied in the middle of the assembly. 
The equation of motion of the system is written as a set of second-order ODEs, i.e., 
\begin{equation}
{\bf M}\ddot{\bf v}(t)+{\bf C}\dot{\bf v}(t)+{\bf K}{\bf v}(t) = {\bf b}(t)+{\bf f}({\bf v}), \label{eq:shell}
\end{equation}
where ${\bf v}$ denotes the nodal displacement vector that contains three translational and two rotaitonal DOFs per node. The nonlinear force ${\bf f}({\bf v})$ is computed based on relative nodal displacements at the debonded area. The force is assumed to be computed by a penalty method, i.e., 
\begin{equation}
f_{i}({\bf v})=-k_p{\rm max}\left\{
0,g(x_u,x_\ell)\right\}
\end{equation}
where $k_p$ denotes a pre-defined penalty parameter and $g$ represents a {\it gap} function defined here as 
\begin{equation}
 g(x_u,x_\ell)\triangleq x_u-x_\ell
 \label{equation:gap}
\end{equation}
where $x_u$ and $x_\ell$ denote the displacements of upper and lower nodes on the surfaces of a specific location of the debonded areas, which are perpendicular to the debonded surfaces. 
\subsection{Galerkin projection based on DMD and constraint modes}\label{subsec:num2:gp}
With FE models {\color{black}{containing multiple contact DOFs}}, it is time consuming to evaluate the contact nonlinearity because contact forces need to be evaluated in physical coordinate system. Therefore, in order to handle the contact problems efficiently, we propose that the DMD modes be used with the so called {\it constraint modes}, which can be obtained by solving a set of static problems where unit displacements are applied to active DOFs that are to be kept in the ROM, with all the other DOFs being fixed~\cite{Craig1968}. 
These modes have been used in the context of {\it Component Mode Synthesis} to apply displacement compatibility conditions to connect multiple bodies. This time, we use them to apply unilateral contact force to the body in the reduced-order space. 
Namely, the dynamics of the entire system is projected onto the subspace spanned by the DMD modes and the constraint modes. Equation (\ref{math:th:-4}) is re-written as follows
\begin{equation}
{\bf x}(t)\approx \sum^{n_m}_{i=1}\eta_i(t)\bm{\varphi}_i + \sum^{n_a}_{i=1}x_i(t)\bm{\psi}_i
\label{numerical_example2:eq2}
\end{equation}
where $\bm{\varphi}_i$ denotes the $i$th DMD mode, $n_m$ denotes the number of DMD modes to be kept in the ROM, $\bm{\psi}_i$ denotes the $i$th constraint mode, $n_a$ denotes the number of constraint modes, which is equal to the number of active DOFs to be kept in the ROM, $\eta_i(t)$ and $x_i(t)$ are the modal coordinates corresponding to the DMD modes and the constraint modes, respectively. 
The modal coordinates $x_i(t)$ behave like physical displacements in the reduced order space. 
This enables ones to handle contact nonlinearity in the reduced order space without expanding the dynamics of the modal coordinates back to the full-order system at every time instant because the DOFs that are subject to the contact nonlinearity are accessible in the reduced-order space. 
\subsection{Effects of debonding on the vibration modes}
\begin{figure}[tb]
\centering
\includegraphics[bb=0 0 1322 1428,width=8.5cm]{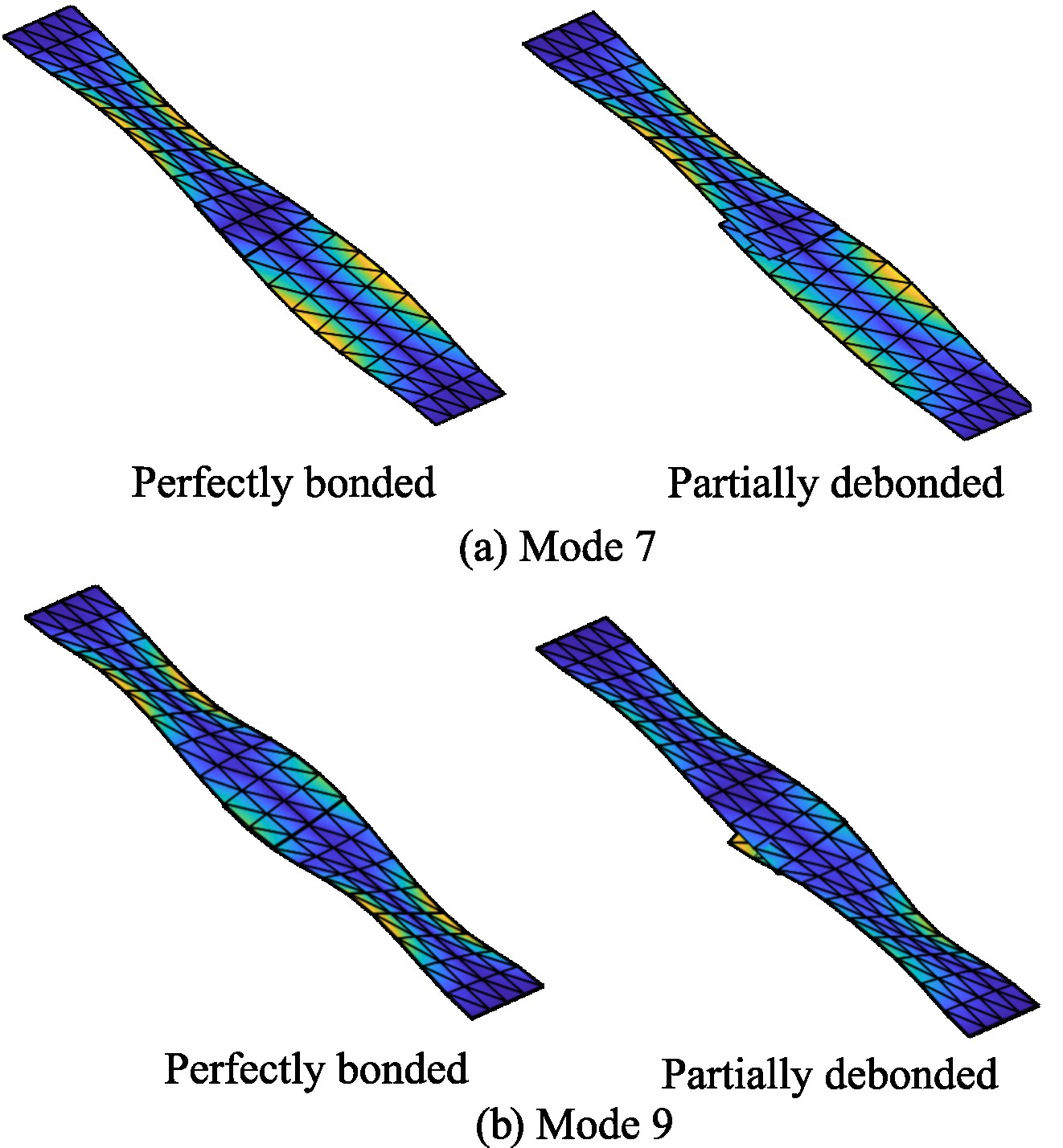}
\caption{Linear normal modes of the assembly}\label{numerical_example2:fig2}
\end{figure}

\begin{table}[tb]
\centering
\caption{Natural frequencies for (A) perfectly bonded and (B) partially debonded assemblies.}\label{numerical_example2:tab1}
\begin{tabular}{cccc}\hline
Mode Number & (A) (Hz) & (B) (Hz) & Deviation (\%)\\ \hline
1&	19.01&	18.93&	-0.40\\
2&	57.02&	56.75&	-0.48\\
3&	114.74&	113.28&	-1.27\\
4&	129.60&	126.87&	-2.11\\
5&	195.78&	194.11&	-0.85\\
6&	311.31&	304.49&	-2.19\\
7&	338.92&	317.08&	-6.45\\
8&	449.00&	444.39&	-1.03\\
9&	522.41&	478.19&	-8.46\\ \hline
\end{tabular}
\end{table}
To illustrate the effects of debonding on the dynamic characteristics of the assembly, linear modal analyses have been conducted on the perfectly bonded shell assembly, and that with partial debonding. Note that with the partial debonding, the contact nonlinearity at the debonded areas has not been enforced, i.e., the system is linear. Also the external force has been ignored for the modal analyses. 
Table \ref{numerical_example2:tab1} shows the first nine natural frequencies for both cases. As can be seen, natural frequencies all decrease because of the debonding, which makes sense because debonding results in softening the stiffness of the system. The effect of debonding is small for the first six and the eighth modes, considering that the deviations in the natural frequencies are less than three percent. On the other hand, the natural frequencies for the seventh and the ninth modes are greatly affected by the debonding considering their deviations are greater than five percent. 
To further examine these variations, the mode shapes corresponding to the seventh and the ninth modes are shown in \figref{numerical_example2:fig2} for both perfectly bonded and partially debonded cases. 
As can be seen in \figref{numerical_example2:fig2}(a), the mode seven shows the twisting shape, whereas the ninth mode also shows the twisting, but higher-order mode shape, as shown in \figref{numerical_example2:fig2}(b).
As can be seen in the mode shapes of the assembly with partial debonding for both modes, a gap at the debonded area appears because the amount of twist for the upper shell is smaller than that for the lower shell at the debonded area. This is because the stiffness with respect to the twisting motion decreased due to the debonding. 

\subsection{Dynamic mode decomposition}
\begin{figure}[tb]
\centering
\includegraphics[scale=1]{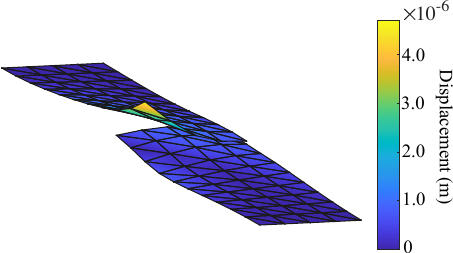}
\caption{Initial deformation of the plate assembly}\label{numerical_example:fig8}
\end{figure}
To obtain DMD modes, snapshots of the displacements need to be obtained. As in the previous numerical example, the snapshots are obtained by solving initial value problems. This time, instead of applying impulsive forcing, initial displacements were applied to the assembly and the resulting free responses were computed. The initial displacements were obtained by applying static forces to open-up debonded areas and solving the corresponding static problem for the deformation of the plate assembly. The shape of the applied initial deformation is shown in \figref{numerical_example:fig8}. This initial condition is employed to ensure the localized motion of the debonded areas is captured in the snapshot. 
With the initial condition, time integration has been conducted by Newmark-$\beta$ method for $0\leqslant t \leqslant 1.0$s. The time integration has been conducted with a fixed time increment of $\Delta t=3.0519\times10^{-5}$s, which means that the dynamic response of frequency up to 16.384~{\rm kHz} are expected to be captured, which is enough for computing DMDs with much low frequencies. 
\begin{figure}[tb]
\centering
\subfigure[Linear case allowing interpenetration of the surfaces]{\includegraphics[scale=1]{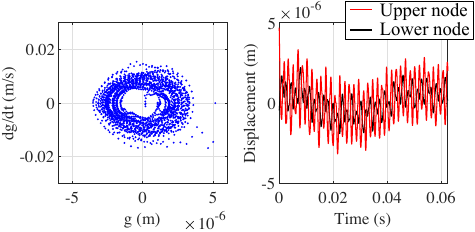}}
\subfigure[Nonlinear case with contact force applied by the penalty method]{\includegraphics[scale=1]{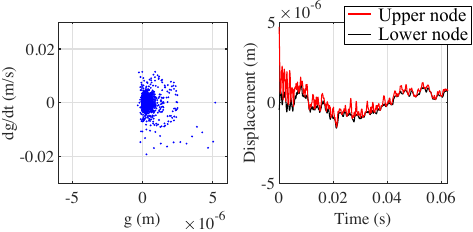}}
\caption{Phase portraits of the gap functions and time histories of the corresponding nodes for $0\leqslant t\leqslant 0.0625$s}\label{numerical_example:fig7}
\end{figure}

The obtained time histories and the corresponding phase portraits of the gap functions are shown in \figref{numerical_example:fig7}, where the gap function is defined as Eq.~\eqref{equation:gap}. 
Figure \ref{numerical_example:fig7}(a) shows a linear response where no contact force is imposed on the potentially contacting nodes even when they are supposed to be in contact. As can be seen in the time histories of the nodes, both nodes oscillate with high frequency with lower node's amplitude being smaller than that for the upper node. The corresponding gap function on the left shows an elliptical orbit with shrinking major and minor axes as time proceeds. 
This is a typical behavior of damped linear systems. 
Figure \ref{numerical_example:fig7}(b) shows the response of the nodes with PWL contact force being enforced. 
As can be seen, the upper node stays above the lower node due to the contact force that occurs when they are in contact, which hinders the interpenetration of the nodes. The trajectory of the gap function in the phase portrait is asymmetric with respect to $g=0$, because gap function is forced to be $g\geqslant 0$ with slight penetration due to the penalty coefficent of finite magnitude. 

\begin{figure}[tb]
\centering
\includegraphics[scale=1]{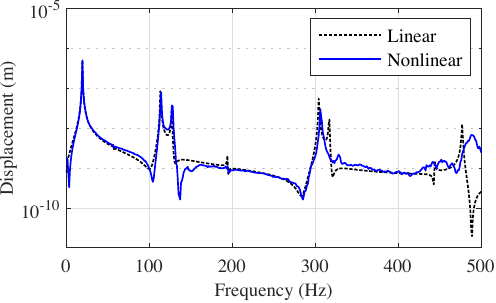}
\caption{Results of FFT of the displacement at the forcing point}\label{numerical_example:fig88}
\end{figure}

The results of the FFT of the displacement at the forcing point for both linear and nonlinear cases are shown in \figref{numerical_example:fig88}. From these results, we can see that the peak frequencies below 200~Hz do not change much regardless of the state of the debonded areas. This agrees with the fact that natural frequencies of the first five modes did not greatly deviate due to the debonding, as shown in \tabref{numerical_example2:tab1}. 
Also, we can see that the peak frequency at 304~Hz for linear case increases to 306~Hz for the nonlinear case, which produced 0.6\% deviation. 
The peak frequency at 317~Hz for the linear case, which corresponds to the seventh mode, increases up to 327~Hz for the nonlinear case, with its amplitude decreased by 13\%. It is because, for the seventh mode, the interpenetration of the debonded areas is so significant without the nonlinear boundary condition that the enforcement of the nonlinear boundary condition at the debonded areas greatly influences the motion of the debonded areas, which results in the suppression of the excitation of such modes. 
\begin{figure}[tb]
\subfigure[Linear (no contact)]{\includegraphics{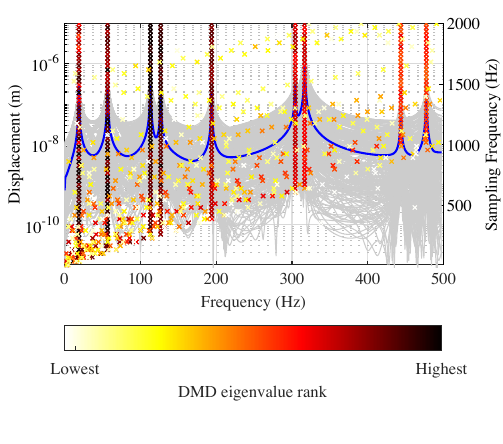}}
\subfigure[Nonlinear (contact at debonded area)]{\includegraphics{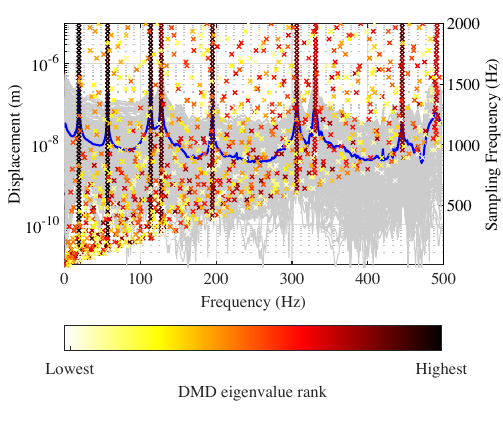}}
\caption{Pseudo-stability diagram of the DMD eigenvalues, along with the average Fourier spectrum of the displacement of all DOFs.  {\color{lightgray}{\bf ---}}: FFT spectra of all DOFs, {\color{blue}{\bf ---}}: average spectrum, ${\bf \times}$: DMD eigenvalues.}\label{numerical_example:plates_pstab}
\end{figure}

\begin{table}
\centering
\caption{DMD frequencies extracted from the free response of the shell assembly by initial displacement}\label{table:shell:dmd}
\begin{tabular}{cccc}\\ \hline
No. & Linear (Hz) & Nonlinear (Hz) & Deviation (\%) \\ \hline
1 &   19.09 &  18.95 &  -0.74 \\
2 &   56.73 &  56.88 &   0.27 \\
3 &  113.25 & 113.64 &   0.34 \\ 
4 &  126.82 & 127.56 &   0.58 \\
5 &  194.05 & 195.07 &   0.52 \\
6 &  304.22 & 306.28 &   0.67 \\
7 &  316.69 & 331.19 &   4.57 \\
8 &  443.53 & 445.28 &   0.39 \\
9 &  477.09 & 491.20 &   2.95 \\
\hline
\end{tabular}
\end{table}

With the computed snapshot, the DMD based on the Algorithm~\ref{alg1} has been applied to both linear and nonlinear cases for increasing sampling frequency to not only capture enough dominant dynamic modes but also to see their physical significance. 
The results of DMD are shown in \figref{numerical_example:plates_pstab}. Again the DMD spectra are shown in conjunction with the FFT spectra of all DOFs and their averages. As can be seen in \figref{numerical_example:plates_pstab}(a), the DMD spectra agree with the peak frequencies in the average FFT spectrum. Moreover, as discussed in \secref{section:example1}, they coincide with the system's natural frequencies when the DMD is applied to the linear response. 
Figure \ref{numerical_example:plates_pstab}(b) shows the DMD spectra computed from the nonlinear case where the contact force at the debonded areas was taken into account. As seen, the dominant DMD frequencies coincide again with the peak frequencies obtained from the FFT. This time, however, they are not necessarily natural frequencies because the response comes from the underlying piecewise linear system. 
In addition, many DMD eigenvalues appear between significant peaks with low rank among the ones computed. This again shows the characteristics of the DMD spectra for nonlinear responses.
Table~\ref{table:shell:dmd} shows the frequencies obtained from the DMD eigenvalues, which are the first nine dominant ones ordered based on the corresponding frequencies. 
As can be seen, most frequencies are not different from each other between the linear and the nonlinear cases, except the mode seven, whose mode shape features localized vibration near the debonded area. 
\begin{figure}[tb]
\centering
\includegraphics[bb=0 0 1004 708,width=8.5cm]{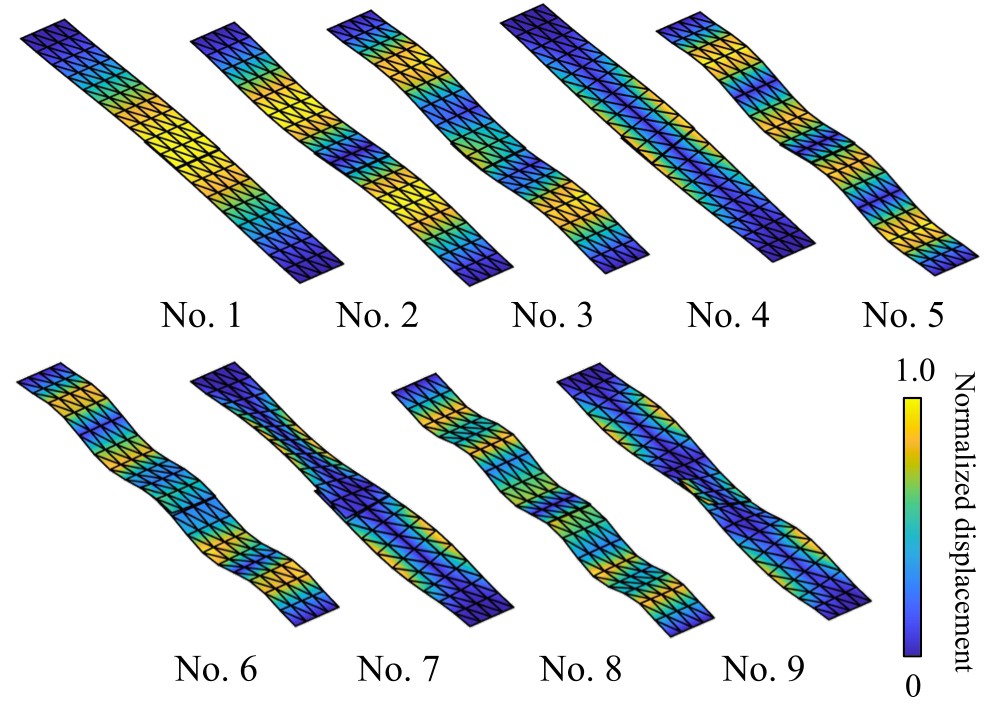}
\caption{Mode shapes of DMD modes ordered based on frequencies of the corresponding DMD eigenvalues}
\label{shell:mode_shape_dmd}
\end{figure}
\begin{figure}[tb]
\centering
\includegraphics[bb=0 0 1004 708,width=8.5cm]{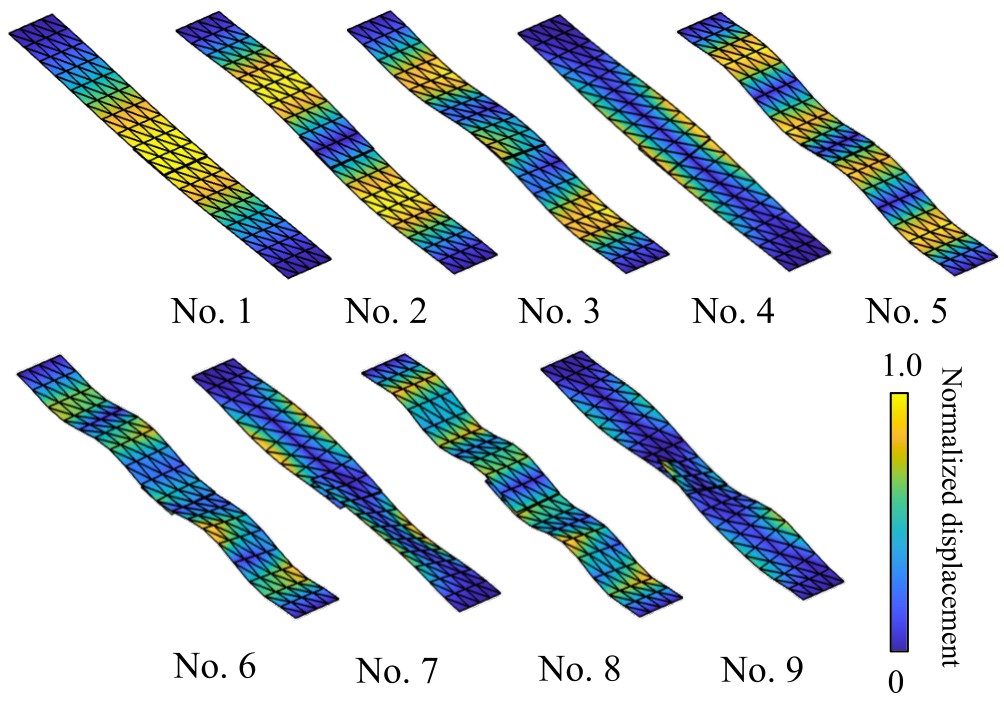}
\caption{Mode shapes of linear normal modes ordered based on frequencies}
\label{shell:mode_shape_lnm}
\end{figure}
\begin{figure}[tb]
\centering
\includegraphics[bb=0 0 1004 708,width=8.5cm]{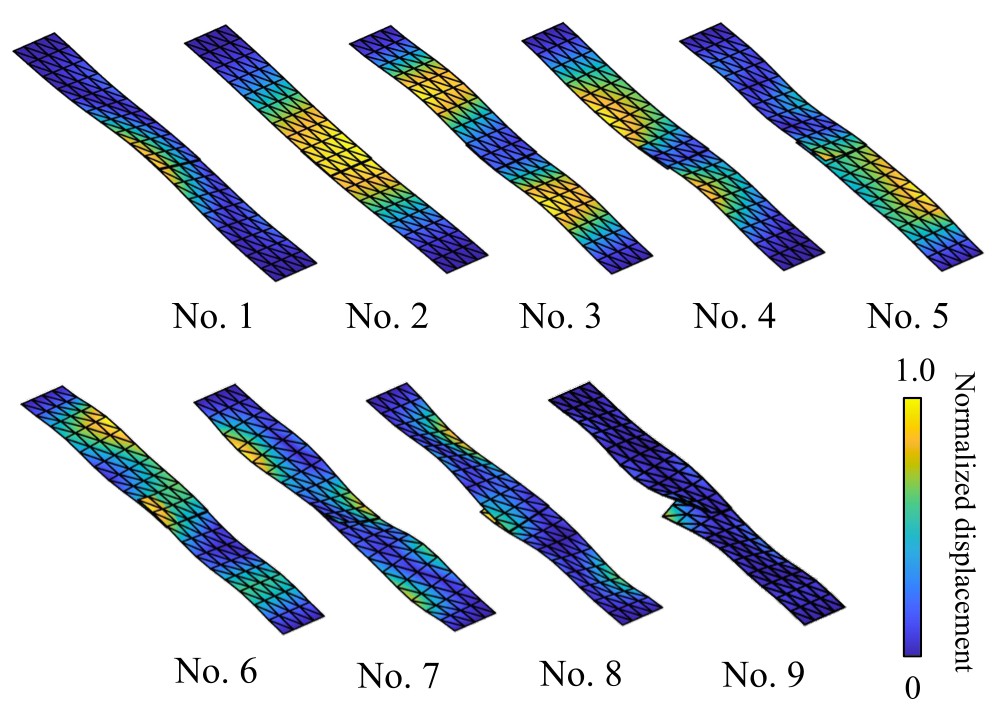}
\caption{Mode shapes of POD modes ordered based on the magnitude of the corresponding singular values $\sigma_i$}
\label{shell:mode_shape_pod}
\end{figure}

%
\subsection{Model order reduction by DMD}
With the obtained DMD modes, Galerkin projection has been applied to Eq.~\eqref{eq:shell}, based on Eq.~\eqref{numerical_example2:eq2}. 
Mode shapes used for the projection are shown in \figref{shell:mode_shape_dmd}. 
In addition, for the sake of comparison, the projection has been applied using LNMs and POD modes. The corresponding shapes of those modes are shown in Figs.~\ref{shell:mode_shape_lnm} and \ref{shell:mode_shape_pod}. 
As can be seen, the mode shapes of DMD modes resemble those of LNMs, with slight differences especially in the vicinity of debonded regions. 
On the other hand, the POD modes do not necessarily resemble neither the LNMs nor the DMDs. Instead, it emphasizes the localized motion at the debonded area. 
For this numerical example, 19 constraint modes corresponding to 18DOFs of the nodes on the debonded surfaces and a DOF corresponding to the excitation point were used for all cases, i.e., $n_a=19$. 

\begin{figure}[tb]
\centering
\includegraphics[scale=1]{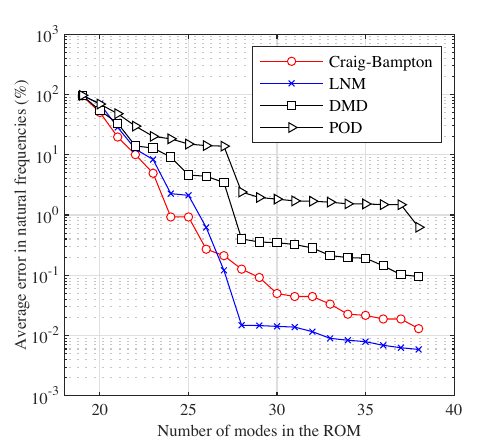}
\caption{Comparison of errors in the natural frequencies for the ROMs compared}
\label{shell:error_in_natural_frequencies}
\end{figure}

First, the obtained ROMs are compared in terms of the errors in their natural frequencies. Namely, the eigenvalues of the ROMs were computed by varying the number of modes in the ROMs. Note that since $n_a=19$ for all cases, the number of modes in Eq.~\eqref{numerical_example2:eq2}, or $n_d$ is varied to see its effects on the ROM accuracy. 
Also, ROMs with the fixed-interface normal modes were considered, which is widely-used accurate model-based MOR and is known to be the Craig-Bampton (CB) method. 
The results are shown in \figref{shell:error_in_natural_frequencies}. As can be seen, LNM converges the fastest with 28 modes. The ROMs with the POD modes also show convergence with 28 modes, but the average error is the largest among the ones compared. The errors computed with the DMD ROMs become approximately 10 times smaller than those obtained with the POD ROMs. The CB-ROM also gives good results. 
Considering that both CB and LNM ROMs are model-based, i.e., all projection basis vectors are computed directly from the system matrices ${\bf M}$ and ${\bf K}$, their capability to capture the eigenspace of ${\bf M}$ and ${\bf K}$ well is not surprising. 
On the contrary, both POD-ROMs and DMD-ROMs are data-driven, although they are compensated by the model-based constraint modes. Therefore, they are not guaranteed to be capable of representing the eigenspace of the system matrices. However, based on these results, we can see that they can produce accurate eigenvalues of the original systems if we include enough number of modes in the projection basis. 

\subsection{Forced response analysis}
\begin{figure*}[tb]
\centering
\includegraphics[scale=1]{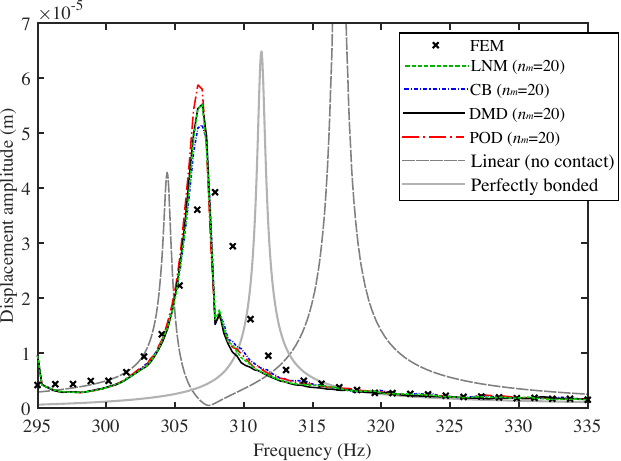}
\caption{Results of forced response analysis}
\label{shell:fr_comparison}
\end{figure*}
Using the FEM model and the ROMs, forced response analysis under harmonic loading has been conducted where $b(t)=A\sin(2\pi f t)$. 
Time integration calculations have been conducted on the equations for $295~{\rm Hz}\leqslant f\leqslant 320~{\rm Hz}$ to excite the sixth and seventh mode until the responses reach the steady-state response. 

The amplitude of the displacement response measured at the excitation point is plotted in \figref{shell:fr_comparison} for FEM, LNM, CB, DMD, and POD ROMs. 
Also, responses of the perfectly bonded assembly and of partially debonded but no contact force is imposed are shown in \figref{shell:fr_comparison}. 
From the plots, we can see that the response of the perfectly bonded assembly shows the largest resonant frequency at 311~Hz. 
On the other hand, partially debonded plate assembly shows the smallest resonant frequency at 304~Hz, if no contact force is applied. 
However, if the contact force that stems from the contact at the debonded region is applied, i.e., if the system is PWL, the responses show the resonance at around 307~Hz. 

Contrary to the expectation from the FFT spectra shown in \figref{numerical_example:plates_pstab} and from the linear response shown in \figref{shell:fr_comparison}, the resonant peak corresponding to seventh mode did not appear. 
Overall, the ROMs predicted the resonance accurately. In terms of the resonance height, CB-ROM produced the best results. Both the LNM and the DMD ROMs produced accurate but slightly larger height than that predicted by the CB-ROM. The resonance height computed by the POD-ROM is the highest and the least accurate. 
\begin{figure}[tb]
\centering
\includegraphics{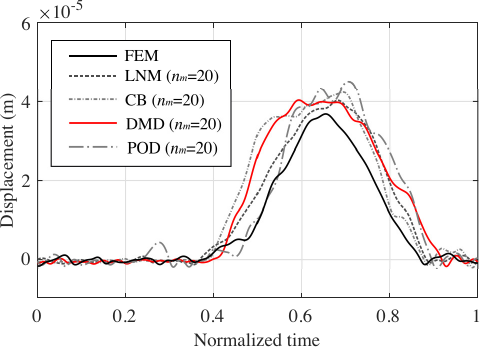}
\caption{Time histories of gap functions at the resonance (307~Hz)}
\label{shell:gap_comparison}
\end{figure}

Figure~\ref{shell:gap_comparison} shows the time histories of the gap functions at the resonance computed by FEM, LNM, CB, DMD, and POD-ROMs. As can be seen, all ROMs accurately predicted the response in comparison with that computed by FEM, but produced slightly larger amplitude than that computed by the FEM. The responses comptued by DMD and LNM produced smooth response that resemble the response computed by the FEM. The responses computed by POD and CB-ROMs contains high-frequency components that do not exist in the one obtained by FEM. 

\begin{figure}[tb]
\centering
\subfigure[Excitation point]{\includegraphics{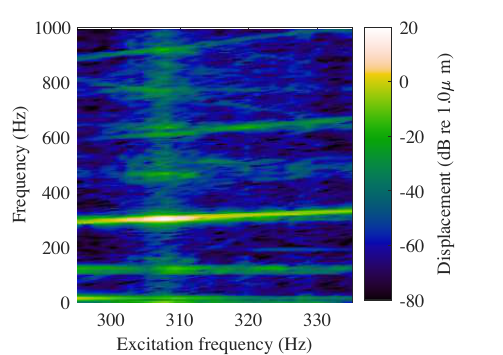}}
\subfigure[Gap function]{\includegraphics{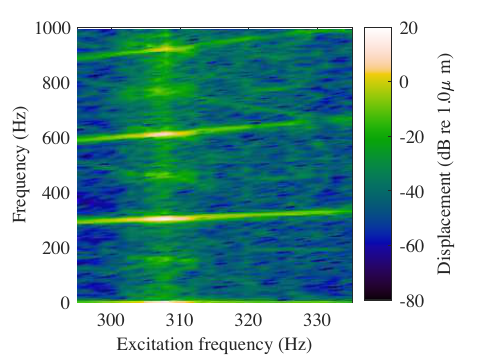}}
\caption{Spectrogram of displacement at the excitation point and gap function computed with FEM}\label{shell:spectrogram_fem}
\end{figure}
\begin{figure}[tb]
\centering
\subfigure[Excitation point]{\includegraphics{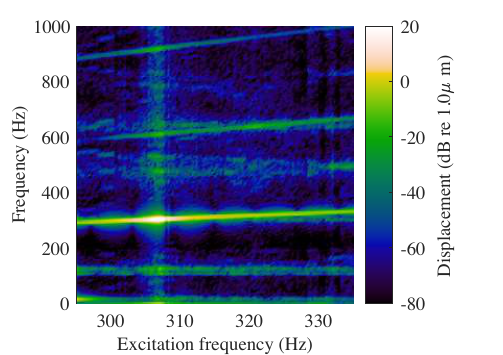}}
\subfigure[Gap function]{\includegraphics{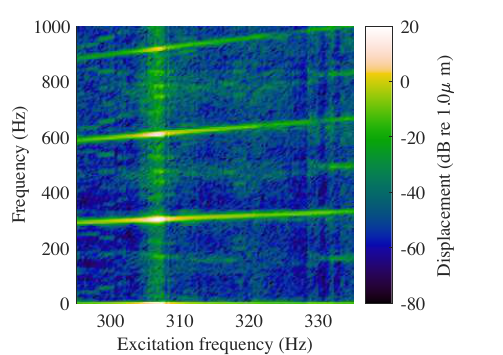}}
\caption{Spectrogram of displacement at the excitation point and gap function computed with DMD-ROM}\label{shell:spectrogram_dmd}
\end{figure}

Furthermore, to visualize the capability of the ROM to represent high frequency components contained in the nonlinear response especially with the DMD modes, spectrogram of the displacement at the excitation point and the gap function was computed by taking FFT on the responses and shown in \figref{shell:spectrogram_fem} and \figref{shell:spectrogram_dmd}. The color represents displacement shown in decibel where 1$\mu$m equals 0dB. 
First, the response at the excitation point is discussed. 
The one predicted by the FEM is shown in Figure~\ref{shell:spectrogram_fem}(a). It shows a strong line corresponding to the first-order of the excitation frequency. The second and the third order lines are also observed, though they are much weaker than the first-order line. Also, there are horizontal lines showing relatively strong response at around 20~Hz, and 120~Hz, which correspond to the first and the fourth modes, respectively, as shown in \tabref{table:shell:dmd}. Note that unlike the order lines, they are not dependent on the excitation frequency. 
When the excitation frequency is at around 307~Hz, there is a resonance, as it was shown in \figref{shell:fr_comparison}. Overall, Figure~\ref{shell:spectrogram_fem}(a) reveals that many linear normal modes and motions with frequencies of integer multiples of excitation frequency occur. 
Figure~\ref{shell:spectrogram_dmd}(a) shows the spectrogram of the displacement evaluated at the excitation point, predicted by the DMD-ROM. Overall, it agrees quite well with the one predicted by the FEM. We can see the strong first-order line and weaker second and third order lines in \figref{shell:spectrogram_dmd}(b). Moreover, the horizontal lines corresponding to the frequencies of the first and the fourth modes can also be observed. The ROM also captures the nonlinear nature of the response at the resonance that contains many linear normal mode components and motions corresponding to the integer-multiples of the excitation frequency. 

Second, the response of the gap function is discussed. 
The one predicted by the FEM, evaluated at a location on the debonded region is shown in \figref{shell:spectrogram_fem}(b). This time, we can see strong, frequency-dependent, first, second, and third order lines in the response. Unlike the response at the excitation point, frequency components corresponding to the linear normal modes are observed only when the excitation frequency is in the vicinity of the resonant frequency. \figref{shell:spectrogram_dmd}(b) shows the one computed by the DMD-ROM. As can be seen, it agrees quite well with the one computed by the FEM. From these, we can see that the response at the debonded area is strongly nonlinear that entails motions with integer-multiples of the excitation frequencies. This is caused by repetitive opening and closing of debonded surfaces. 
%

\begin{table*}[tb]
\centering
\caption{Comparison of computational time required to complete a set of forced response calculations at the resonance (second)}\label{table:cputime}
{\color{black}{
\begin{tabular}{lccccc}\hline \hline
&FEM & LNM & CB & DMD & POD\\ \hline
(A) Snapshot extraction    & N/A & N/A & N/A & $6.7\times10^3$ &$6.7\times10^3$\\ 
(B) DMD (POD) calculations & N/A & N/A & N/A & $5.9\times10^1$ &$1.2\times10^{-1}$\\ 
(C) Galerkin projection    & N/A &$7.2\times10^{-3}$ & $1.2\times10^{-2}$ & $5.8\times10^{-3}$ & $4.8\times10^{-3}$\\ 
(D) Forced response calculation (per frequency) & $7.6\times10^3$ & $2.7\times10^1$ & $1.8\times10^1$ & $2.2\times10^1$ & $2.7\times10^1$\\ 
Total (A)+(B)+(C)+(D)$\times n_f$                          & $9.7\times10^5$ & $3.5\times10^3$ & $2.3\times10^3$ & $9.6\times10^3$ & $1.0\times10^4$\\ \hline
\end{tabular}
}}
\end{table*}

Lastly, the benefit of the proposed MOR methodology in terms of the computational time is discussed. The forced response calculations were conducted by the Newmark $\beta$ time integration method with the ROMs and the FEM for $f=307$~Hz and $0\leqslant t\leqslant 0.208$s, where the number of time steps was fixed to 65,536 for all cases, resulting in the fixed time step size of $\Delta t=3.17\times10^{-6}\mbox{s}$. 
Also, $n_m=20$ for all the ROMs considered. 
The forced response calculations were repeated five times with the same conditions and the average computational time for each ROM {\color {black}{per frequency}} was computed. 
The measurements were conducted on a workstation with the CPU of AMD Threadripper PRO 3975WX (4.2~GHz) and 128GB of RAM. During the measurements, the application of the numerical simulation used only a single core of the CPU for all cases. 
The results are shown in \tabref{table:cputime} {\color {black}{where CPU time for each computational step is shown for FEM, LNM, CB, DMD and POD. The total CPU time shown in the bottom row is shown with the assumption that the CPU time for forced response calculation per frequency is constant throughout the frequency range. }}
As can be seen in the table, although the number of modes was chosen to be the same, there are slight differences between the computational times for the ROMs. {\color{black}{
In particular, LNM and CB show smaller computational times than DMD and POD. This is because they are classical model-based ROM that does not require the computation of the snapshots. On the other hand, both DMD and POD take longer CPU time than the model-based ROMs but they are comparable. This means that with the proposed method, we are able to generate ROMs with a reasonable amount of time and accuracy. }}
Furthermore, we can see that the computations with the ROMs were much faster than that with the FEM. Indeed, the computations with DMD-ROM were approximately 100 times faster than those with the FEM model. 
From these, we can see the efficiency of using the ROMs in comparison with the full-order FEM model when predicting the forced response of the system with PWL nonlinearity.  

\section{Conclusion}\label{sec:conclusion}
In this paper, a model order reduction method for piecewise-linear systems based on data-driven dynamic mode decomposition has been proposed. 
The key idea of the proposed method is that the snapshots are obtained based on impulse response or response due to initial deformation of the system, which results in rich spectral contents in the response due to the piecewise-linear nonlinearity. This gives one the capability to capture the complicated dynamics due to PWL nonlinearity, which normally contains many frequency components. 
For the systems with many DOFs with piecewise-linear nonlinearity, it was proposed that DMD modes be used with the constraint modes to be able to handle the nonlinearity efficiently. 
The proposed methodology has been applied to a forced response problems of a cantilevered beam with an elastic stop at the tip, and a bonded shell assembly with partial debonding subjected to harmonic loading. 
The forced response problems were solved with the proposed method, and the results agreed well with the ones obtained by the full-order FEM models. 

%
%
%
%
\section*{Statements and Declarations}

\noindent{\bf Funding}\\
This project was supported in part by Japan Society for the Promotion of Science (JSPS), Grant-in-Aid for Scientific Research(C), grant number JP20K11855. 
\\ \\
\noindent{\bf Competing interests}\\
The authors have no relevant financial or non-financial interests to disclose. 
\\ \\
\noindent{\bf Author contributions}\\
Both authors contributed to the study conception, design and analysis. Simulations were implemented and performed mainly by A. Saito. The first draft of the manuscript was written by A. Saito and M. Tanaka commented on the draft of the manuscript. All authors read and approved the final manuscript. 
\\ \\
\noindent{\bf Data availability}\\
The dataset generated during and/or analysed during the current study are available from the corresponding author upon reasonable request. 

\end{document}